\documentclass{amsart}

\usepackage{amsmath}
\usepackage{amssymb}
\usepackage{mathrsfs}
\usepackage{url}
\usepackage{enumitem}
\usepackage{stmaryrd}
\usepackage[all]{xy}
\usepackage{verbatim}
\usepackage{color}

\newtheorem{Thm}{Theorem}[section]

\newtheorem{Cor}[Thm]{Corollary}

\newtheorem{Prop}[Thm]{Proposition}
\newtheorem{Qn}[Thm]{Question} 
\theoremstyle{definition}

\newtheorem{Rmk}[Thm]{Remark}

\numberwithin{equation}{Thm}

\newcommand{\CC}{\mathbb{C}}

\newcommand{\FF}{\mathbb{F}}

\newcommand{\QQ}{\mathbb{Q}}

\newcommand{\ZZ}{\mathbb{Z}}
\newcommand{\PP}{\mathbb{P}}
\newcommand{\A}{\mathbb{A}}

\newcommand{\ml}[2]{\begin{multline}\label{#1}#2 \end{multline}}
\newcommand{\ga}[2]{\begin{gather}\label{#1}#2 \end{gather}}

\date{\today}

\author{H\'el\`ene Esnault}
\address{Freie Universit\"at Berlin, Arnimallee 3, 14195, Berlin, Germany}
\email{esnault@math.fu-berlin.de}
\author{Daqing Wan}
\address{Department of
Mathematics, University of California, Irvine, CA 92697, USA}
\email{dwan@math.uci.edu}


\title{Divisibility of Frobenius Eigenvalues on $\ell$-adic cohomology}

\begin{document}

\begin{abstract}
For a projective variety defined over a finite field  with  $q$ elements, it is shown that  as algebraic integers, 
the eigenvalues of the geometric Frobenius acting 
on $\ell$-adic cohomology
have higher  than known $q$-divisibility 
beyond the middle dimension. This sharpens both Deligne's integrality theorem  \cite[Corollary 5.5.3]{SGA7} and 
 the  cohomological divisibility theorem  \cite[Theorem~4.1]{EK05}.  Similar lower bound is proved 
 for the Hodge level for a complex projective variety beyond the middle dimension, improving earlier results in this direction. 
The affine version of our results for the compactly supported cohomology is still open in general.

\end{abstract}

\maketitle

\section{Introduction}

In this   small note, we continue to explore the interplay between the Ax-Katz theorem \cite{Ka71} for point count $q$-divisibility and 
Deligne's integrality theorem \cite{SGA7} for varieties defined  over a finite field $\FF_q$ with  $q$ elements. 
This connection was studied earlier in   \cite{Es03}  for projective complete intersections and  \cite{EK05} 
for both projective and affine varieties. Building on these earlier works,  we show that for a projective variety, there is an  improved $q$-divisilibity  of 
the eigenvalues of the geometric Frobenius acting on $\ell$-adic cohomology beyond the middle dimension.  The proof carries over to $\CC$  to yield a lower bound  for the Hodge level for 
projective varieties over $\CC$,  improving earlier results in \cite{Es90}, \cite{ENS92}, \cite{EW03}. 

We first describe the affine case and then return to the projective case. 
Let $X$ be an affine variety in $\A^n$ over $\FF_q$, defined by an ideal  spanned by $r$ polynomials 
$$f_1, \ldots, f_r \in \FF_q[x_1,\cdots, x_n]$$
of positive degree $d_1, \ldots, d_r$.  
To avoid triviality, we assume that $r$ and $n$ are positive integers,  
thus  in particular $\dim(X) \leq n-1$. 
For each non-negative integer $j$, we define another non-negative integer 
$$ \mu_j(n; d_1,..., d_r) = j + \max( 0, \lceil \frac{n-j- \sum_{i=1}^r d_i}{\max_{i=1}^r d_i}\rceil ) \geq j.$$
This is an increasing function in both $j$ and $n$, that is
\ml{}{\mu_{j+1}(n; d_1,\ldots, d_r) \geq \mu_j(n; d_1,\ldots, d_r), \\ \mu_{j}(n+1; d_1,\ldots, d_r)  \geq \mu_j(n; d_1,\ldots, d_r).\notag}
It is clearly a decreasing function in each $d_i$. 
The Ax-Katz theorem \cite{Ka71} says that the number $\# X(\FF_q)$ of $\FF_q$-rational points on $X$ is divisible by 
$q^{\mu_0(n;  d_1,\ldots, d_r)}$ for all $q$. This yields a  non-trivial information as soon as $n > \sum_{i=1}^r d_i$.
Equivalently,  as $\mu_0(n; d_1,\ldots, d_r)\le n$, $\# (\A^n \setminus X)(\FF_q)$ is divisible by 
$q^{\mu_0(n; d_1,\ldots, d_r)}$.  
 The bound is achieved,  that is the result is sharp,
 so this $q$-divisibility result is best possible in general.  It has a natural interpretation in terms 
of zeta functions. 

Recall that the zeta function of $X$ is defined by the formal power series
$$Z(X, T) = \exp( \sum_{k=1}^{\infty} \frac{\# X(\FF_{q^k})}{k}T^k)   \in  1+ T \ZZ[[T]]\subset \QQ[[T]],  $$
which by  Dwork's rationality theorem \cite{Dw60} lies in $\QQ(T)$.
This implies that the reciprocal zeros and poles of $Z(X,T)$ are algebraic integers. 
From the $q$-adic radius of convergence 
for the logarithmic derivative of the zeta function, 
one deduces that the Ax-Katz theorem 
is equivalent to the statement that all reciprocal zeros and poles of $Z(X,T)$
are divisible by $q^{\mu_0(n; d_1,\ldots, d_r)}$ as algebraic integers. 

 By Grothendieck's trace formula \cite{Gr65}, the zeta function has a cohomological interpretation
$$Z(X,T)= \prod_{i\ge 0} P_{2i+1}(T)/\prod_{i\ge 0} P_{2i}(T), $$
where $P_i(T)={\rm det} (I-TF|H^i_c(X))$
and by a (standard) abuse of notation
$ H^i_c(X):=H^i_c(X\otimes \bar \FF_q, \QQ_\ell). $
Here $\ell$ is a  prime  different from the characteristic $p$ 
of $\FF_q$, $H^i_c$ denotes  the $i$-th $\ell$-adic compactly supported cohomology, $F$ is the geometric Frobenius. 
By  the Weil conjectures, that is Deligne's purity theorem \cite{De74}, $q^{\mu}$-divisibility  as algebraic integers  of the  eigenvalues of the geometric Frobenius  
is equivalent to $q^{\kappa \mu}$-divisibility of  
$\# X(\FF_{q^\kappa})$ for all $\kappa \geq 1$ if $X$ is smooth proper. In general, from the divisibility of the number of points, one cannot immediately deduce the divisibility of the eigenvalues of  the geometric Frobenius as there might be some cancellation in the expression of the zeta function as a rational function.  
 It is then natural to ask what can be said about the $q$-divisibility 
on  $H^j_c(X)$.  This question was first studied in \cite{Es03} for projective complete intersections  
and then in \cite{EK05} in the general case. 
We state their general results  for both affine and projective varieties. The affine result is
\begin{Prop}\label{Prop1}
Let $X$ be an affine variety in $\A^n$   defined by  $r$ polynomials $f_1,\ldots, f_r \in \FF_q[x_1,\ldots, x_n]$ of 
positive 
degrees $d_1,\ldots, d_r$  with $n, r\ge 1$.  
 \begin{itemize}
 \item[(i)]  \cite[Theorem 2.1 and Theorem 2.2]{EK05} For all integers $i$, the
 eigenvalues of the geometric Frobenius acting on $H_c^{i}(X)$ and $H_c^{i}(\A^n\setminus X)$ 
 are divisible by $q^{\mu_0(n; d_1,\ldots, d_r)}$ as algebraic integers. 
 \item[(ii)]  \cite[Theorem 2.3]{EK05} For all integers $j\geq 0$, the
 eigenvalues of the geometric Frobenius acting on $H_c^{n-1+j}(X)$ and $H_c^{n+j}(\A^n\setminus X)$ 
 are divisible by $q^{\mu_j(n; d_1,\ldots, d_r)}$ as algebraic integers.  \end{itemize}

\end{Prop}

The projective result is

\begin{Prop}\label{Prop2}\cite[Theorem 4.1]{EK05}  
Let $Y\subset \PP^n$ be a projective variety defined by $r$ homogeneous polynomials 
$f_1,\ldots, f_r \in \FF_q[x_0, x_1,\ldots, x_n]$of positive degrees $d_1,\ldots, d_r$ with $n, r\ge 1$.

\begin{itemize} 
\item[(i)] For all integers $i$, the eigenvalues of the geometric Frobenius acting on $H_c^{i}(\PP^n\setminus Y)$ 
are divisible by $q^{\mu_0(n+1; d_1,\ldots, d_r)}$ as algebraic integers. 

\item[(ii)] For all integers $j\geq 0$, the eigenvalues of the geometric Frobenius acting 
$H_c^{n+1+j}(\PP^n \setminus Y)$ are divisible by $q^{\mu_j(n+1; d_1,\ldots, d_r)}$ as algebraic integers.\end{itemize}

\end{Prop}
In the case when $Y$ is a complete intersection in $\PP^n$, part (i)  of Proposition \ref{Prop2}  was first proved in \cite{Es03}. 
This part (i)  gives a cohomological  strengthening  of the Ax-Katz theorem which   in the projective case  is the statement that the eigenvalues  of the geometric Frobenius 
on  $H_c^{i}(\PP^n\setminus Y)$ 
are divisible by $q^{\mu_0(n+1; d_1,..., d_r)}$ as algebraic integers for all $i$.  
 The proof uses
the Ax-Katz theorem and thus does not reprove it.  
Part (ii)  of Proposition \ref{Prop2}  further 
says that for $j\geq 0$, the eigenvalues of the geometric Frobenius on $H_c^{n+1+j}(\PP^n \setminus Y)$ are divisible by 
$q^{\mu_j(n+1; d_1,..., d_r)}$ as algebraic integers.  This is better than part (i), since 
$\mu_j(n+1; d_1,..., d_r) \geq \mu_0(n+1; d_1,..., d_r)$.

 Our main result is the following projective theorem which improves part (ii) of Proposition \ref{Prop2}.

\begin{Thm}\label{Thm1.3}  Let $Y\subset \PP^n$ be a projective variety defined by $r$ homogeneous polynomials 
$f_1,\ldots, f_r \in \FF_q[x_0, x_1,\ldots, x_n]$of positive degrees $d_1,\ldots, d_r$ with $n, r\ge 1$. We denote by $H^i_{\rm prim}(Y)=H^i(Y)/H^i(\PP^n)$ the primitve cohomology. Then the eigenvalues of the geometric Frobenius acting on 

\begin{itemize} 
\item[(i)] $H_{\rm prim} ^{\dim(Y)+j}(Y)$  for $0\leq j \leq \dim(Y)$,
\item[(ii)] $H_c^{\dim(Y)+1+j}(\PP^n \setminus Y)$  for $0\leq j \leq \dim(Y)+1$
\end{itemize}
are divisible by $q^{\mu_j(n+1; d_1,\ldots, d_r)}$ as algebraic integers.

\end{Thm}
 This improves part (ii) of Proposition \ref{Prop2} 
because now the higher divisibility starts at ${\rm dim}(Y)+1$, earlier than $n+1$. Even at the same cohomological degree $n+1+j$, 
Theorem \ref{Thm1.3} would give higher divisibility than Proposition \ref{Prop2} does, since 
$$H_c^{n+1+j}(\PP^n \setminus Y) = H_c^{\dim(Y)+1+ (j+n-\dim(Y))}(\PP^n \setminus Y).$$
In the non-trivial case where $\max_i d_i >1$, one checks that for $0\leq j\leq \dim(Y)$, 
$$\frac{\dim(Y)+j}{2} \geq \mu_j(n+1; d_1, \ldots, d_r),$$
so we can replace the primitive cohomology in Theorem \ref{Thm1.3} (i) by the ordinary 
cohomology. We state this remark  as a corollary.

\begin{Cor} In Theorem~\ref{Thm1.3} , if  $\max_i d_i >1$, then 
for all $j\geq 0$, the eigenvalues of the geometric  Frobenius acting on 
$H^{\dim(Y)+j}(Y)$ 
are divisible by $q^{\mu_j(n+1; d_1,\ldots, d_r)}$ as algebraic integers. 
\end{Cor}

Since $\mu_j(n+1; d_1,\ldots, d_r) \geq j$, this corollary also improves Deligne's integrality theorem  
\cite[Corollary 5.5.3]{SGA7} which says that the eigenvalues of the geometric  Frobenius acting on 
$H^{\dim(Y)+j}(Y)$ 
are divisible by $q^j$ as algebraic integers. 

In this note, we have improved the projective Proposition \ref{Prop2} beyond middle cohomological 
dimension. 
It would be interesting to know if the affine Proposition \ref{Prop1} can be similarly improved beyond middle 
cohomological dimension. We state this as an open 
problem. 

\begin{Qn}\label{Qn}
Let $X$ be an affine variety in $\A^n$   defined by  $f_1,\ldots, f_r \in \FF_q[x_1,\ldots, x_n]$ of 
positive 
degrees $d_1,\ldots, d_r$  with $n, r\ge 1$. Is it true that the
 eigenvalues of the geometric Frobenius acting on 
 \begin{itemize}
 \item[(i)]  $H_c^{\dim(X)+j}(X)$ for $j\ge 0$,
 \item[(ii)] $H_c^{\dim(X) +1+j}(\A^n\setminus X)$ for $0\le j \le {\rm dim}(X)+1$ \end{itemize}
are divisible by $q^{\mu_j(n; d_1,\ldots, d_r)}$ as algebraic integers?  
\end{Qn}
\begin{Rmk} \label{rmk:zhang} 
The projective Theorem \ref{Thm1.3} can already be used to treat one non-trivial case. 
Let $Y$ be the projective variety in 
$\PP^n$ defined the vanishing of the homogeneous polynomials 
$$g_i(x_0, x_1, \ldots, x_n): =x_0^{d_i}  f_i(\frac{x_1}{x_0}, \dots, \frac{x_n}{x_0}),  \ 1\leq i \leq r.$$
Let $Y_{\infty}$ be the hyperplane section $Y\cap \{x_0=0\}$ at infinity,  so 
$$\PP^n \setminus Y = (\A^n \setminus X) \sqcup (\PP^{n-1} \setminus Y_{\infty}).$$
By the excision sequence

\ga{}{\ldots \to  H^j_c(\A^n\setminus X)\to H^j_c(\PP^n\setminus Y)\to H^j_c(\PP^{n-1} \setminus Y_{\infty}) \to \ldots \notag}the affine case for   $\A^n \setminus  X$  is reduced to the projective case for both  $\PP^n \setminus Y$ and 
 $\PP^{n-1} \setminus Y_{\infty}$
 in case  $\dim(Y)=\dim(X)$. 
Thus, Question~\ref{Qn} has a positive answer if $\dim(Y)=\dim(X)$. 

In general, our inductive proof shows that Question~\ref{Qn} has a positive answer if for a general hyperplane 
$A$ in $\A^n$, the affine Gysin map 
$$H^{i-2}_c(A \setminus  A\cap X)(-1)\xrightarrow{\rm Gysin} H^i_c(\A^n \setminus X) $$
was surjective for $i-1>{\rm dim}(X)$. 
Its projective version is   true. 
 If $X$ is a local complete intersection in $\A^n$,  the affine 
Gysin lemma
 follows from the Weak Lefschetz Theorem for perverse sheaves  \cite[A.5]{Del93}
as observed by Dingxin Zhang.
Thus, Question~\ref{Qn} also has a positive answer if $X$ 
is a local complete intersection in $\A^n$. 
\end{Rmk}

 An immediate application of  the above remark  is a cohomological  strengthening 
 of the polar result in \cite[Theorem~1.2b]{Wa02} for affine complete intersections.  A non-zero number $\alpha$ is called a reciprocal zero (resp. reciprocal pole) of the rational function $Z(X,T)$ if $1/\alpha$ is a zero (resp. pole) of 
$Z(X,T)$. 
In the case that $X$  is an affine complete intersection in $\A^n$, then $H_c^i(X)=0$ for all $i< \dim(X)$ and 
the cohomological formula for the zeta function reduces to 
$$Z(X,T)^{(-1)^{\dim(X)-1}}= \prod_{j=0}^{\dim(X)} \det(I -  TF | H_c^{\dim(X)+j}(X))^{(-1)^{j}}.$$
This implies that the reciprocal poles of $Z(X,T)^{(-1)^{\dim(X)-1}}$ are among the Frobenius eigenvalues 
on $H_c^{\dim(X)+j}(X)$ for odd $j$, and hence are divisible as algebraic integers by 
$$q^{\min_{j\geq 0}\mu_{2j+1}(n; d_1,..., d_r)} =  q^{\mu_{1}(n; d_1,..., d_r)}.$$
This recovers the polar result in \cite[Theorem~1.2b]{Wa02} according to which 
if $X$ is an affine complete intersection in $\A^n$,  all reciprocal poles of $Z(X,T)^{(-1)^{\dim(X)-1}}$
are divisible by $q^{\mu_1(n; d_1,..., d_r)}$ as algebraic integers,   where we note that 
$$\mu_1(n; d_1, \cdots, d_n) = 1 + \mu_0(n-1; d_1, \cdots, d_n)\geq \mu_0(n; d_1, \cdots, d_n).$$

The proof of Theorem \ref{Thm1.3}  does not  contain any new geometric or cohomological idea. 
It depends on earlier results in \cite{Es03} and \cite{EK05}. 
For the proof of Theorem~\ref{Thm1.3} (ii),  we mimic the method in \cite[Section~2]{Es03}.
The whole point is to have the right numbers which allow an inductive argument and apply the result 
in [EK05] as our starting point (the case $j=0$).
In Section~\ref{sec:Hodge}, we mention  the Hodge theoretic analogues of our divisibility 
result. 

\medskip
{\bf Acknowledgements}. We are pleased to thank Nick Katz for helpful discussions on this topic.   We are also grateful 
 to  the referee  for his/her comments which led to an improved exposition  and to Dingxin Zhang for a precise reading, see Remark~\ref{rmk:zhang}.

\section{Proofs }
\subsection{Proof of Theorem~\ref{Thm1.3}}
We prove Theorem~\ref{Thm1.3} (ii).  We argue by induction on ${\rm dim}(Y)$. If ${\rm dim}(Y)=-1$ that is $Y=\emptyset$, then $r>n$ thus $\mu_0(n+1; d_1,..., d_r)=0$ and there is nothing to prove.
 We  assume that $\dim(Y)\geq 0$. If $j=0$, 
the theorem is  already true by  Proposition~\ref{Prop2} (i). We  assume that $\dim(Y)\geq 0$ and $j\geq 1$. 
If $n=1$, then  ${\rm dim}(Y)$  must be zero and $j=1$. In this case, Theorem~\ref{Thm1.3} (ii) is true as well 
since $\mu_1(2, d_1,\cdots, d_r)=1$ and $H_c^2(\PP^1\setminus Y)$ has the unique Frobenius eigenvalue $q$. 
Thus, we assume that $\dim(Y)\geq 0$, $j\geq 1$ and $n\geq 2$ below.   
  
We copy the argument of   \cite[Section~2]{Es03}.  By replacing $\FF_q$ by a finite extension, which preserves the divisibility, there is by \cite[Theorem~2.1]{Es03} a linear hyperplane $\iota: A\hookrightarrow  \PP^n$ such that the Gysin homomorphism
$$H^{i-2}(A, \iota^*\mathcal F)(-1)\xrightarrow{\rm Gysin}  H^i_A(\PP^n, \mathcal F)$$ is a Frobenius equivariant isomorphism, where $\mathcal F=j_!\mathbb Q_\ell, \ j: \PP^n\setminus Y\to \PP^n$. On the other hand, one has an exact sequence 
\ga{}{ H^{i-1}(Y\setminus Y\cap A)\to H^i(\PP^n\setminus A, \mathcal F)=H^i(\PP^n\setminus A, Y\setminus Y\cap A)\to H^i(\PP^n\setminus  A) \notag}
so Artin vanishing on the left and right terms  implies 
$$
 H^i(\PP^n\setminus A, \mathcal F)=0 \ {\rm for} \ i-1> {\rm dim}(Y).$$
Thus the composite 
\ml{}{ H^{i-2}_c(A \setminus  A\cap Y)(-1)= H^{i-2}(A, \iota^*\mathcal F)(-1) \xrightarrow{\rm Gysin} H^i_A(\PP^n, \mathcal F)\notag\\ \xrightarrow{\rm excision} H^i(\PP^n, \mathcal F) =H^i_c(\PP^n \setminus Y) }
is surjective
and Frobenius equivariant for $i-1>{\rm dim}(Y)$.
 We now assume this inequality and set $i={\rm dim}(Y)+1+j$  so $j\ge 1$.
 Since ${\rm dim}(Y)\geq 0$ and $A$ is in general position, we have $\dim(A\cap Y) =\dim(Y)-1$. 
By our induction hypothesis, the theorem is true for the lower dimensional projective variety $A\cap Y$ in $A =\PP^{n-1}$ 
defined by the vanishing of $r$ homogeneous polynomials of degrees $d_1, \cdots, d_r$ with $(n-1), r\geq 1$. 
It follows that for $0\leq j-1\leq \dim(A\cap Y)+1$, 
the eigenvalues of the geometric Frobenius on $H_c^{\dim(A\cap Y)+1+(j-1)}(A\setminus A\cap Y)(-1)$ 
are divisible by 
$$ q^{1+\mu_{j-1}(n; d_1,\ldots,  d_r)}=q^{\mu_{j}(n+1; d_1,\ldots, d_r)}$$ as algebraic integers. 
So  for $1\leq j\leq \dim(A\cap Y)+2 =\dim(Y)+1$, 
the divisibility  on $H_c^{\dim(Y)+1+j}(\PP^n\setminus  Y)$ 
is  also  by $q^{\mu_{j}(n+1; d_1,\ldots, d_r)}$.  This finishes the proof of Theorem~\ref{Thm1.3} (ii).  

 To prove Theorem~\ref{Thm1.3} (i),   one applies the excision sequence which gives the isomorphism of Frobenius modules
$$H_{\rm prim} ^{i}(Y): = H^i(Y)/H^i(\PP^n) \cong H_c^{i+1}(\PP^n-Y), \ 0\leq i \leq 2\dim(Y)\leq 2(n-1).$$
This implies that for $0\leq j \leq \dim(Y)$, Frobenius eigenvalues on
$$H_{\rm prim} ^{\dim(Y)+j}(Y)= H^{\dim(Y)+j}(Y)/H^{\dim(Y)+j}(\PP^n) \cong H_c^{\dim(Y)+1+j}(\PP^n-Y)$$
are divisible by $q^{\mu_{j}(n+1, d_1,..., d_r)}$ as algebraic integers. The theorem is proved.

\section{Hodge level} \label{sec:Hodge}
The divisibility as algebraic integers of the Frobenius eigenvalues on $\ell$-adic cohomology suggests  by a vast generalization of Tate conjecture
a 
similar divisibility of the corresponding motive  by the Tate motive. The motivic divisibility in turn implies
a similar lower bound for the Hodge level, which we sketch in this section.

For a non-empty separated finite type scheme $X$ over $\CC$, 
the compactly supported cohomology group $H_c^i(X)$ has a mixed Hodge structure  \cite{De71}  with the decreasing 
Hodge filtration ${\rm Fil}$: 
$$H_c^i(X) = {\rm Fil}^0 H_c^i(X)\supseteq {\rm Fil}^1 H_c^i(X)\supseteq {\rm Fil}^2 H_c^i(X) \supseteq \cdots.$$
The Hodge level of $H_c^i(X)$ is the largest integer $\kappa $ such that 
$${\rm Fil}^\kappa H_c^i(X) = H_c^i(X).$$

\begin{Thm}\label{H}  Let $Y\subset \PP^n$ be a projective variety defined over $\CC$ by $r$ homogeneous polynomials 
of positive degrees $d_1,\cdots, d_r$ with $n, r\ge 1$. We denote by $H^i_{\rm prim}(Y)=H^i(Y)/H^i(\PP^n)$ the primitve cohomology. The 
Hodge level of 

\begin{itemize} 
\item[(i)] $H_{\rm prim} ^{\dim(Y)+j}(Y)$  for $0\leq j \leq \dim(Y)$,
\item[(ii)] $H_c^{\dim(Y)+1+j}(\PP^n \setminus Y)$  for $0\leq j \leq \dim(Y)+1$
\end{itemize}
is  at least $\mu_j(n+1; d_1,\ldots, d_r)$.
\end{Thm}

The proof of this theorem is the same as that of Theorem~\ref{Thm1.3}. 
Indeed one has 
Artin's vanishing theorem, excision, purity and base change. One does the induction in the same way. 
It reduces  to the case $j=0$.  
For $j=0$, the Hodge level lower bound is proved in \cite{ENS92} in the 
projective case.
 The latter extends earlier results in \cite{DD90} for projective hypersurfaces and  in \cite{Es90} 
for projective complete intersections. 
 
 \begin{Rmk} In the complete intersection case, Theorem~\ref{H} (ii) was first proved in \cite[Theorem 2.3]{EW03} 
 using a different geometric argument. The complete intersection condition comes from a vanishing theorem 
 for Zariski sheaf cohomology in \cite{Es90}. This suggests that it might be possible to refine 
 the geometric approach in \cite{DD90}, \cite{Es90}, \cite{EW03} to prove Theorem~\ref{H} as well for any projective variety $Y$ in $\PP^n$. 
 
In the affine case, we expect  a similar improvement and we state this as an open problem. 

 \begin{Qn}\label{Q1}
Let $X$ be an affine variety in $\A^n$ defined over $\CC$ by  $r$ polynomials  of 
positive 
degrees $d_1,\ldots, d_r$ with $n, r \ge 1$. Is it true that 
the Hodge level of 
 \begin{itemize}
 \item[(i)]  $H_c^{\dim(X)+j}(X)$ for $j\ge 0$,
 \item[(ii)] $H_c^{\dim(X) +1+j}(\A^n\setminus X)$ for $0\le j \le {\rm dim}(X)+1$ \end{itemize}
is at least  $\mu_j(n; d_1,\ldots, d_r)$?  
\end{Qn}

\end{Rmk}


\begin{thebibliography}{9999}

\bibitem[De71]{De71} \textsc{P. Deligne},  Th\'eorie de Hodge, II,  {\it Publ. Math. IHES} {\bf 40} (1971), 5--58.


\bibitem[De74]{De74} \textsc{P. Deligne},  La conjecture de Weil : I, 
{\it Publ Math IHES}, {\bf 43} (1974), 273-307. 




\bibitem[Del93]{Del93}  \textsc{P. Deligne}, Appendix to: Affine cohomological transforms, perversity, and monodromy, by  \textsc{N. Katz}, {\it J.  Amer. Math. Soc.} {\bf 1} 6 (1993), 218--222. 

\bibitem[DD90]{DD90}  \textsc{P. Deligne, A. Dimca},  Filtrations de Hodge et par l'ordre du pôle pour les hypersurfaces singulières,
{\it Ann. Sci.  \'Ecole Norm. Sup.},  {\bf 23} (1990), 645-656. 

\bibitem[Dw60]{Dw60} \textsc{B. Dwork}, On the rationality of the zeta function of an algebraic variety, 
{\it Amer. J. Math.}, {\bf 82} (1960), 631-648. 

\bibitem[Es90]{Es90} \textsc{H. Esnault}, 
Hodge type of subvarieties of $\PP^n$  of small degrees, {\it  Math. Ann.},   {\bf 288 }(1990), 549--551 

\bibitem[ENS92]{ENS92} \textsc{H. Esnault, M. Nori, V. Srinivas},  Hodge type of projective varieties of low degree, {\it Math. Ann.}, {\bf 293} (1992), 1--6.


\bibitem[Es03]{Es03} \textsc{H. Esnault}, Eigenvalues of Frobenius acting on the $\ell$-adic cohomology of complete intersections of
low degree, 
{\it C. R. Math. Acad. Sci. Paris},  {\bf 337} (2003), 317-320.




\bibitem[EK05]{EK05} \textsc{H. Esnault, N. Katz}, 
Cohomological divisibility and point count divisibility, 
{\it Compositio Math.}, {\bf 141}(2005),  93-100. 


\bibitem[EW03]{EW03} \textsc{H. Esnault,  D. Wan}, Hodge type of the exotic cohomology of complete intersections, 
{\it C. R. Math. Acad. Sci. Paris},  {\bf 336} (2003), 153-157.


\bibitem[Gr65]{Gr65} \textsc{A. Grothendieck}, Formule de Lefschetz et rationalit\'e des fonctions L, S\'eminaire Bourbaki, 
Vol 9, Expos\'e {\bf 279} (1964/65), 41-55. 


\bibitem[Ka71]{Ka71} \textsc{N. Katz}, On a theorem of Ax, 
{\it Amer. J. Math.},  {\bf 93} (1971), 485-499.


\bibitem[Wa02]{Wa02} \textsc{D. Wan}, Poles of zeta functions of complete intersections, 
{\it Chinese Ann. Math.}, {\bf  21B} (2002), 187-200.


\bibitem[SGA4]{SGA4} \textsc{ M. Artin}, Th\'eor\`eme de changement de base par un morphisme propre, et applications, in S\'eminaire de G\'eom\'etrie Alg\'ebrique {\bf 4}, Tome 3, Expos\'e XVI, 
{\it Lecture Notes in Mathematics} {\bf 305}, Springer Verlag (1973), 206--249.
 
 \bibitem[SGA4.5]{SGA4.5} \textsc{P. Deligne}, Th\'eor\`emes de finitude en cohomologie $\ell$-adique,  in S\'eminaire de G\'eom\'etrie Alg\'ebrique {\bf 4 1/2}, {\it Lecture Notes in Mathematics} {\bf 569}, Springer Verlag (1977), 233--251.
 
 \bibitem[SGA7]{SGA7} \textsc{P. Deligne}, Th\'eor\`eme d’int\'egralit\'e, Appendix to N. Katz, Le niveau de la 
cohomologie des intersections compl\`etes, Expos\'e XXI, in SGA 7, 
{\it Lecture Notes in Mathematics}, {\bf Vol. 340}, Springer Verlag (1973), 363-400.


\end{thebibliography}
\end{document}